\documentclass{amsart}

\usepackage{amsmath,amssymb,amsthm}
\usepackage{graphicx,tikz}

\newtheorem*{proposition}{Proposition}

\theoremstyle{definition}

\theoremstyle{remark}

\begin{document}

\title[]{Localization of Neumann eigenfunctions\\ near irregular boundaries}

\author[]{Peter W. Jones}
\address{Department of Mathematics, Yale University, New Haven, CT 06511, USA}
\email{jones@math.yale.edu}

\author[]{Stefan Steinerberger}
\address{Department of Mathematics, Yale University, New Haven, CT 06511, USA}
\email{stefan.steinerberger@yale.edu}

\thanks{S.S. is supported by the NSF (DMS-1763179) and the Alfred P. Sloan Foundation.}

\begin{abstract} It has been empirically observed that eigenfunctions of Laplace's equation $-\Delta \phi = \lambda \phi$ with Neumann
boundary conditions sometimes localize near the boundary of the domain if that boundary is rough (say, fractal). This has some nontrivial implications in acoustics that has been put to real-life use
(sound attenuation by noise-protective walls); this short paper describes the mathematical mechanism responsible for this and describes the quantitative strength of the phenomenon for some examples.
\end{abstract}

\maketitle

\vspace{-10pt}

\section{Introduction}
\subsection{Introduction.} This short paper is concerned with an interesting localization phenomenon that has been observed to occur for eigenfunctions
of the Laplace operator with Neumann boundary conditions. Let us fix an open, bounded domain $\Omega \subset \mathbb{R}^n$ (our subsequent argument also works on manifolds) and consider
\begin{align*}
 -\Delta u &= \mu u \qquad \mbox{inside}~\Omega \\
\frac{\partial u}{\partial \eta}&=0 \qquad \quad \mbox{on} ~\partial \Omega.
\end{align*}
It is known that there is a discrete sequence of eigenvalues $0 = \mu_0 < \mu_1 \leq \mu_2 \leq \dots$ for which there is a solution and we will refer to these solutions
as the eigenfunctions of the Neumann Laplacian.
Their localization behavior is a fundamental concept which we do not summarize here (see \cite{agmon, anderson, arnold, fil, fil2, greb, ng2, steinerberger}
and references for examples of classical and recent work on this). It is classical that if $\Omega$ is comprised of essentially two domains connected
by a thin connecting set, then the first Neumann eigenfunction is essentially localized in each cluster; this is valuable in applications \cite{chung, coifman}.
We focus on a specific subtle phenomenon that has been actively studied by B. Sapoval, collaborators and others \cite{fel1, fel2, hab, heilman, ng2, ng, sap1, sapoval, sap2, sap3, sap4}.

\begin{center}
\begin{figure}[h!]
\begin{tikzpicture}[scale=2.5]
\draw [ultra thick] (0,0.3) -- (1,0.3) -- (1,1) -- (0.9, 0.8) -- (0.85, 1.03) -- (0.77, 0.85) -- (0.7, 1.1) -- (0.65, 0.8) -- (0.57, 1) -- (0.5, 0.9) -- (0.45, 0.8) -- (0.4, 1) --(0.37, 0.8) -- (0.32, 0.9) -- (0.3, 0.88) -- (0.27, 1) -- (0.23, 0.78) -- (0.18, 1) -- (0.13, 0.85) -- (0.05, 1.1) -- (0,1) -- (0,0.3);
\draw [ultra thick]  (2, 0.3) -- (3,0.3) -- (3,1.2) -- (2,1.2) -- (2,0.3);
\filldraw [ultra thick]  (2.2, 0.8) -- (2.8, 0.8) -- (2.8, 0.85) -- (2.2, 0.85) -- (2.2, 0.792);
\filldraw [ultra thick]  (2.2, 1) -- (2.8, 1) -- (2.8, 1.05) -- (2.2, 1.05) -- (2.2, 0.99);
\end{tikzpicture}
\caption{Eigenfunctions have been observed \cite{fel1} to localize close and inside the irregular boundary (left) and between slits (right).}
\end{figure}
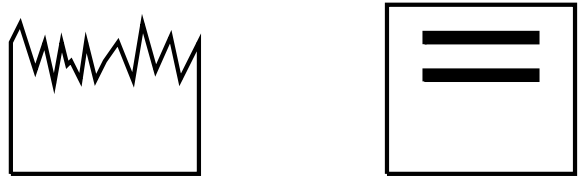
\end{center}

A suitable example, taken from Felix, Asch, Filoche \& Sapoval \cite{fel1}, is shown in Figure 1: for these two domains, one can numerically observe the existence of localized Neumann eigenfunctions close to the irregular boundary (in the 'spikes'
and between the slit, respectively). A survey of Grebenkov \& Nguyen \cite[Section 7.6.]{greb} discusses the phenomenon in greater detail. As is observed in \cite{fel1}, the presence of irregular boundaries seems to cause some of the Neumann
eigenfunctions to localize there even though these irregularities do not necessarily separate a subdomain from the rest of the domain in any obvious way: 
this has some nontrivial implications and is even used in the construction of road noise barriers (\cite{fel1} refers to the Fractal Wall$^{\tiny \mbox{TM}}$, product of Colas Inc., French patent No. 0203404).\\

What complicates matters, again refering to the domain shown in Fig. 1, is that the localization is not as pronounced as in classical examples in mathematical physics where a localized eigenfunction undergoes exponential decay away from the area where it is localized. Here, there is no exponential decay in the rest of the
domain -- moreover, the localized eigenfunctions on the domain with spikes do not necessarily localize within a single 'spike' but in several neighboring spikes at the same time. The cause of the phenomenon has been puzzling for a while: originally it was believed to be a phenomenon mainly found in fractal domains. 
A 2010 \textsc{Notices} articles of Heilman \& Strichartz \cite{heilman} gives non-fractal examples and instead puts an emphasis on the role of symmetries; the recent survey of Grebenkov \& Nguyen \cite{greb}
remarks 'symmetry is neither sufficient nor necessary for localization'.  Heilman \& Strichartz also remark that 'We
do not have to go very high up in the spectrum' [to find these localized eigenfunctions].

\section{A Proof of Localization}
\subsection{Idea.} There is a fairly simple explanation for the phenomenon: the irregular boundary 'traps' Brownian motion (which is reflected) for an extended period of time
and this leads to the creation of localized eigenfunctions. 

\begin{center}
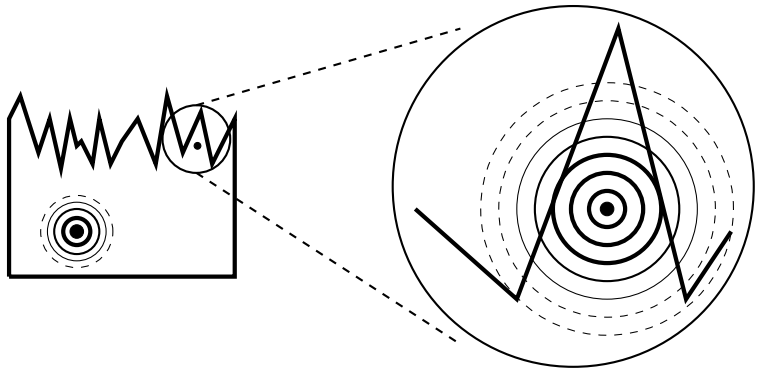
\begin{figure}[h!]
\begin{tikzpicture}[scale=3]
\draw [ultra thick] (0,0.3) -- (1,0.3) -- (1,1) -- (0.9, 0.8) -- (0.85, 1.03) -- (0.77, 0.85) -- (0.7, 1.1) -- (0.65, 0.8) -- (0.57, 1) -- (0.5, 0.9) -- (0.45, 0.8) -- (0.4, 1) --(0.37, 0.8) -- (0.32, 0.9) -- (0.3, 0.88) -- (0.27, 1) -- (0.23, 0.78) -- (0.18, 1) -- (0.13, 0.85) -- (0.05, 1.1) -- (0,1) -- (0,0.3);
\draw [thick]  (0.83, 0.91) circle (0.15cm);
\draw [thick, dashed] (0.83,1.06) -- (2, 1.4);
\draw [thick, dashed] (0.83,1.01-0.25) -- (2, 0);
\draw [thick] (2.5,0.7) circle (0.8cm);
\draw [ultra thick] (1.8, 0.6) -- (2.25, 0.2) -- (2.7, 1.4) -- (3, 0.2) -- (3.2, 0.5);
\filldraw (0.3, 0.5) circle (0.03cm);
\draw [ thick] (0.3, 0.5) circle (0.1cm);
\draw [ultra thick] (0.3, 0.5) circle (0.06cm);
\draw [] (0.3, 0.5) circle (0.13cm);
\draw [dashed] (0.3, 0.5) circle (0.16cm);
\filldraw (0.835, 0.88) circle (0.015cm);
\filldraw (2.65, 0.6) circle (0.03cm);
\draw [ultra thick] (2.65, 0.6) circle (0.08cm);
\draw [ultra thick] (2.65, 0.6) circle (0.16cm);
\draw [ultra thick] (2.65, 0.6) circle (0.24cm);
\draw [ thick] (2.65, 0.6) circle (0.32cm);
\draw [] (2.65, 0.6) circle (0.4cm);
\draw [dashed] (2.65, 0.6) circle (0.48cm);
\draw [dashed] (2.65, 0.6) circle (0.56cm);
\end{tikzpicture}
\caption{Heat Kernel inside the domain (left) and close to the irregulary boundary (right): the Heat kernel gets reflected back many times leading
to a higher concentration close to the point.}
\end{figure}
\end{center}
 We will now explain the underlying idea and apply it to the two domains in Figure 1. It is not hard to use the idea to derive explicit statements in various other settings.
Let $p(t,x,y)$ denote the Neumann heat kernel
$$ p(t,x,y) = \sum_{k=0}^{\infty}{e^{-\mu_k t} \phi_{k}(x) \phi_{k}(y) }.$$
We recall that $p(t,x,y) \geq 0$, that for all $x \in \Omega$ and all $t>0$
$$ \int_{\Omega}{ p(t,x,y) dy} = 1$$
and that $p(t,x,y) = p(t,y,x)$. In free Euclidean space, $p(t,x,\cdot)$ is merely a Gaussian centered at $x$ and localized at spatial scale $\sim \sqrt{t}$. The probabilistic interpretation is that of $p(t,x,\cdot)$ as describing the distribution of a Brownian motion particle started at $x$, running for time $t>0$ and being reflected upon impact on the boundary. For very irregular boundaries, there might be some difficulty in defining the notion of a reflected Brownian motion, however, we only work in settings where Neumann eigenfunctions of the Laplacian are defined.
The key identity that underlies the phenomenon is 
\begin{align*}
\int_{\Omega}{ p(t,x,y)^2 dy} &= \int_{\Omega}{ \left( \sum_{k=0}^{\infty}{e^{-\mu_k t} \phi_{k}(x) \phi_{k}(y) } \right)^2 dy}  \\
&= \int_{\Omega}{\sum_{k, \ell=0}^{\infty}{e^{-\mu_k t} e^{-\mu_{\ell} t} \phi_{k}(x) \phi_{k}(y) \phi_{\ell}(x) \phi_{\ell}(y)}  dy} \\
&= \sum_{k=0}^{\infty}{ e^{-2 \mu_k t} \phi_k(x)^2}.
\end{align*}
This identity can be used in both directions: if $p(t,x, \cdot)$ is strongly concentrated around $x$ (which happens if the boundary $\partial \Omega$ has a trapping effect on Brownian motion), then this implies that at least some
of the Neumann eigenfunctions have to be large in $x$. Moreover, these eigenfunctions have to correspond to a relatively small eigenvalue, large values arising for large eigenvalues are dampened the exponential decay and play no role. Conversely, if $p(t,x,\cdot)$ is spread out, then none of the low-frequency eigenfunctions can be concentrated in $x$ beyond a certain degree. The parameter $t>0$ is arbitrary and has the effect of fixing a natural length scale $\sim \sqrt{t}$ and an associated frequency range of eigenvalues that play a role.

\section{Scaling and some examples} The purpose of this Section is to point out the natural scales of these objects in the case of $\Omega \subset \mathbb{R}^2$ since all the examples in this paper are two-dimensional. However, we hasten to emphasize that the higher-dimensional case (or the manifold setting) can be dealt with in exactly the same way.
 Ignoring (small) contributions from possibly fractal boundaries
(see \cite{berry, lap}), we expect $\mu_k \sim ck$ on two-dimensional domains with an implicit constant $c > 0$ depending only on the measure $|\Omega|$ of the domain. We also
observe that in the case of all eigenfunctions being relatively flat, i.e. $\phi_k(x)^2 \sim 1$, we expect that
$$ \sum_{k=1}^{\infty}{e^{-\mu_k t} \phi_k(x)^2}  \sim \sum_{k=1}^{\infty}{e^{-\mu_k t}} \sim  \sum_{k=1}^{\infty}{e^{- ck t}} = 1 + \frac{1}{e^{ct} -1} \sim \frac{1}{ct}.$$
Naturally, this chain of equivalences can be made in a variety of possible ways under a variety of different assumptions.
We will now compute the left-hand side for domains $\Omega \subset \mathbb{R}^2$ that behave like the figure on the left-hand side of Fig. 1 (a square with thin cones attached).
We first compute it for points $x$ inside
$\Omega$ that are far away from the boundary $d(x, \partial \Omega) \gtrsim \sqrt{t}$. Points far away from the boundary have a short-time heat kernel asymptotic dictated by Varadhan's Lemma \cite{vara1}.
This suggests, in two dimensions, that
$$ p(t,x,y) \sim  \frac{1}{(4 \pi t)^{}} e^{-\frac{\|x-y\|^2}{4t}} \qquad \mbox{and thus} \qquad \int_{\Omega}{ p(t,x,y)^2 dy} \sim \frac{1}{8 \pi t}$$
This is very much in correspondence with what we would expect under the assumption that $\phi_k(x)^2 \sim 1$ for all $k$. However,
the situation inside a cone is different.

\begin{center}
\begin{figure}[h!]
\begin{tikzpicture}[scale=3]
\draw [ultra thick] (0.7, 0.7) -- (0,0) -- (1,0);
 \draw [thick, domain=0:45] plot ({0.7*cos(\x)}, {0.7*sin(\x)});
\node at (0.4, -0.15) {$\sqrt{t}$};
\node at (0.3, 0.5) {$\sqrt{t}$};
\node at (0.25,0.1) {$\alpha$};
\filldraw (0.15, 0.08) circle (0.02cm);
\node at (1.3, 0.3) {$\implies$};
\draw [ultra thick] (2.7, 0.7) -- (2.35,0.35) -- (3,0.35);
\draw [thick] (2.35, 0.35) circle (0.35cm);
\draw [thick] (2.35, 0.35) -- (2.35,0.35) -- (2.35, 0.7) -- (2.35,0.35) -- (2, 0.35) -- (2.35,0.35) -- (2, 0.7) -- (2.35,0.35) -- (2,0) -- (2.35,0.35) -- (2.35, 0) -- (2.35, 0.35) -- (2.7, 0) -- (2.35,0.35);
\end{tikzpicture}
\caption{Reflection trick.}
\end{figure}
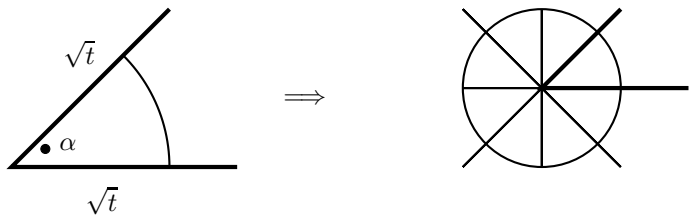
\end{center}

We use the reflection trick to compute the integral for $x$ at distance $\leq 0.01 \sqrt{t}$ from the apex of a cone with opening angle $\alpha$ and sides of length at least $\sim \sqrt{t}$. We obtain, in order of magnitude,
$$  \int_{\Omega}{ p(t,x,y)^2 dy} \sim \frac{2\pi}{\alpha}\frac{1}{8 \pi t} = \frac{1}{\alpha}\frac{1}{ 4  t}.$$
This, however, is much larger than the corresponding quantity in free-space if $\alpha$ is small and implies
$$  \sum_{k=0}^{\infty}{ e^{-2 \mu_k t} \phi_k(x)^2} \sim \frac{1}{\alpha} \frac{1}{4t}$$
which forces some of the low-lying eigenfunctions to be large in $x$. We observe that this does \textit{not} imply exponential localization (which
would be false) nor does it imply that there is an eigenfunction that is large close to the apex of exactly one cone (which would also be false): it only
implies that there exists an eigenfunction for which $e^{-\mu_k t}\phi_{\mu}(x)^2 \gtrsim \alpha^{-1}$ which is a factor of $\sim \alpha^{-1}$ larger than one would expect for a function
at scale $\sim 1$. This is only possible for $\mu_k$ small and thus also
explains why these localized objects are observed at rather low frequencies. It also shows that the phenomenon does not operate a different scale
but rather at the same scale with a larger leading constant.

\subsection{Another example: the slit.} We illustrate this fundamental principle in a second case: that of a slit introduced in a two-dimensional domain $\Omega$
(see Fig. 1). Let us assume, for simplicity, that the slits are of size $\sim 1 \times 0.01$ each and that they are distance $\delta$
apart from each other.

Let us pick $x$ to be the point in the center between the two slits and let us fix $t = 1$. We observe that the likelihood of Brownian
motion escaping the regions $R$ between these two slits, a region of area $\sim \delta$, is bounded away from 0. 
\begin{center}
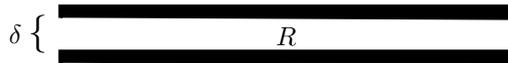
\begin{figure}[h!]
\begin{tikzpicture}[scale=3]
\filldraw [ultra thick] (0, 0) -- (0,0.05) -- (2,0.05) -- (2,0) -- (0,0.01);
\filldraw [ultra thick] (0, 0.2) -- (0,0.25) -- (2,0.25) -- (2,0.2) -- (0,0.21);
\node at (-0.1, 0.13) {\LARGE $\{$};
\node at (-0.2, 0.13) {$\delta$};
\node at (1, 0.12) {$R$};
\end{tikzpicture}
\caption{A thin slit in an otherwise nice domain.}
\end{figure}
\end{center}
This means that
$$ 1 \lesssim \int_{R}{ p(1,x,y) dy} \leq  \delta^{1/2} \left( \int_{R}{ p(1,x,y)^2 dy} \right)^{\frac{1}{2}}$$
and thus
$$  \sum_{k=0}^{\infty}{ e^{-2 \mu_k} \phi_k(x)^2} \gtrsim \frac{1}{\delta}.$$
This determines the scale at which the localization phenomenon increases as the two slits move closer to each other. All of this can
be made precise under some assumption on the domain $\Omega$ (which determines the scaling of $\mu_k$).

\subsection{Non-Localization.}
There is also a converse implication whenever a point $x$ is far from the boundary, thus allowing for an unobstructed diffusion of the
Green's function and implies that
$$  \int_{\Omega}{ p(t,x,y)^2 dy} \sim \frac{1}{t} ~\qquad \quad \mbox{will behave as in free space}.$$
In that case, none of the first few eigenfunctions can be very large in $x$ (unless many more vanish which would imply them localizing somewhere else). This has an interesting implication for domains such as the ones displayed
in Fig. 1: the boundary is actually piecewise-linear comprised of finitely many pieces. This shows that the only point in which diffusion can
truly be trapped in the sense above, is close to the apex of a cone for the domain shown on the left in Fig. 1. Moreover, diffusion cannot be traped at all in the slit domain as soon
as $t \lesssim \delta^2$. For the domain with a slit this reasoning suggests, under some
control on the eigenvalues and assuming an otherwise smooth boundary $\partial \Omega$, that above a certain frequency no
eigenfunction will be particularly localized. Conversely, in the presence of cone structures the argument suggests the existence of infinitely
many eigenfunctions that have substantial size close to the apex.

\subsection{Another Example: Fractals.} Returning to the original question of the behavior of eigenfunctions near the boundary of irregular or possibly even
fractal domains, we see that the main question is whether diffusion gets trapped or whether it spreads freely. It is natural to ask whether
the speed of diffusion (and thus, by the argument above, the strength of the localization of Neumann eigenfunctions) is somehow related to the Hausdorff dimension of the fractal. \\

We consider one explicit example for the sake of clarification by taking a sequence of disks of radius $r_n = 2^{-n}$ that are glued together along line segments of length $2^{-4^{n}}$ (see Fig. 5).

\begin{proposition} Some eigenfunctions exhibit exponential growth in the balls.
\end{proposition}

We observe that the same result is true for eigenfunctions of the Dirichlet Laplacian for rather trivial reasons: the narrow connection
between two adjacent balls essentially separates these domains completely, eigenfunctions localize and thus scale with the inverse volume that decays exponentially - hence exponential growth. Eigenfunctions of the Neumann-Laplacian behave differently, there is no reason why they would be small in a slit connecting two balls.
\begin{center}
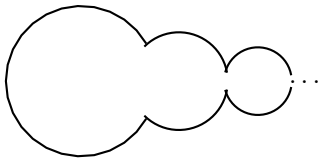
\begin{figure}[h!]
\begin{tikzpicture}
   \draw [thick,domain=30:330] plot ({cos(\x)}, {sin(\x)});
      \draw [thick,domain=10:135] plot ({1.3+0.65*cos(\x)}, {0.65*sin(\x)});
      \draw [thick,domain=180+45:350] plot ({1.3+0.65*cos(\x)}, {0.65*sin(\x)});
      \draw [thick,domain=10:165] plot ({2.35+0.45*cos(\x)}, {0.45*sin(\x)});
      \draw [thick,domain=195:350] plot ({2.35+0.45*cos(\x)}, {0.45*sin(\x)});
      \node at (3,0) {$\dots$};
\end{tikzpicture}
\caption{Balls with exponentially decaying radii glued together along even fast exponentially decaying line segments.}
\end{figure}
\end{center}

\begin{proof} The argument is a slight modification of the previous arguments (but still relying on the same mechanism). Pick $\chi_n$ to be the characteristic function of the $n-$th ball. We apply the Neumann heat equation for time $t=1$ and obtain
$$ e^{\Delta}\chi_n =  \sum_{k=0}^{\infty}{ e^{- \mu_k } \left\langle  \chi_n, \phi_k \right\rangle \phi_k}.$$
It follows from known estimates on the narrow escape problem \cite{singer} that for this choice of parameters $e^{t \Delta} \chi_n \geq \chi_n/2$. The Cauchy-Schwarz inequality implies
$$ \left| \left\langle  \chi_n, \phi_k \right\rangle \right| \leq \|\chi_n\|_{L^2} \| \phi_k\|_{L^2} = \| \chi_n\|_{L^2} \leq c \cdot 2^{-n}.$$
Altogether, plugging in a point $x$ from the $n-$th ball, we obtain
$$ \frac{1}{2} \leq e^{\Delta}\chi_n(x) =   \sum_{k=0}^{\infty}{ e^{- \mu_k } \left\langle  \chi_n, \phi_k \right\rangle \phi_k} \leq
 c \cdot 2^{-n} \sum_{k=0}^{\infty}{ e^{- \mu_k } \left| \phi_k(x) \right|}.$$
 This shows
  $$ \sum_{k=0}^{\infty}{ e^{- \mu_k } \left| \phi_k(x) \right|} \geq  2^{n-1}/c \qquad \mbox{for}~x~\mbox{in the $n$-th ball}.$$
\end{proof}

\subsection{Related results.}
The entire approach is in a similar spirit as work of the second author \cite{steinerberger} giving an interpretation of localization properties
in terms of diffusion. Moreover, the same functional has recently been used by J. Lu and the second author \cite{lu} where a discrete heat
kernel was used to determine the location of localized eigenvectors associated to linear system $Ax = \lambda x$. In particular, 
$$  \int_{\Omega}{ p(t,x,y)^2 dy}  ~ \mbox{can be used to detect localized eigenfunctions.}$$
The paper \cite{lu} also discusses fast linear algebra methods for its computation.

\end{document}